\documentclass[journal]{IEEEtran}
\usepackage{lineno,hyperref}
\usepackage{amsfonts}
\usepackage{amssymb}
\usepackage{amsmath}
\usepackage{graphicx} 
\usepackage{enumerate}
\usepackage{booktabs} 
\usepackage{hypernat} 
\usepackage{stfloats}

\newcommand \T {\mathrm{T}}

\usepackage{cite}
\ifCLASSINFOpdf
\else
\fi
\begin{document}
\title{Distributed Model Predictive Control Under Inexact Primal-Dual Gradient Optimization Based on Contraction Analysis}

\author{Yanxu~Su,
	~Yang~Shi,~\IEEEmembership{~Fellow,~IEEE,}
	and~Changyin~Sun% <-this % stops a space
	\thanks{This work was supported in part by the National Natural Science Foundation of China under Grant 61520106009 and Grant U1713209.}
	\thanks{Y. Su and C. Sun are with the School of Automation, Southeast University, Nanjing 210096, China, and also with the Key Laboratory of Measurement and Control of Complex System of Engineering, Ministry of Education,
		Southeast University, Nanjing 210096, China (e-mail: yanxu.su@seu.edu.cn; cysun@seu.edu.cn).}% <-this % stops a space
	\thanks{Y. Shi is with the Department of Mechanical Engineering, University of Victoria, Victoria, BC V8W 3P6, Canada (e-mail: yshi@uvic.ca).}}% <-this % stops a space



\markboth{Journal of \LaTeX\ Class Files,~Vol.~14, No.~8, August~2015}%
{Shell \MakeLowercase{\textit{et al.}}: Bare Demo of IEEEtran.cls for IEEE Journals}



\maketitle

\begin{abstract}
This paper develops a distributed model predictive control (DMPC) strategy for a class of discrete-time linear systems with consideration of globally coupled constraints. The DMPC under study is based on the dual problem concerning all subsystems, which is solved by means of the primal-dual gradient optimization in a distributed manner using Laplacian consensus. To reduce the computational burden, the constraint tightening method is utilized to provide a capability of premature termination with guaranteeing the convergence of the DMPC optimization. The contraction theory is first adopted in the convergence analysis of the primal-dual gradient optimization under discrete-time updating dynamics towards a nonlinear objective function. Under some reasonable assumptions, the recursive feasibility and stability of the closed-loop system can be established under the inexact solution. A numerical simulation is given to verify the performance of the proposed strategy.
\end{abstract}


% \fi
\IEEEpeerreviewmaketitle

\section{Introduction}
%\IEEEPARstart{T}{his} demo file is intended to serve as a ``starter file''

Model predictive control (MPC) is of tremendous interest in recent years for its broad applications ranging from industrial control systems \cite{qin2003survey}, aerospace \cite{eren2017model}, sensor network control \cite{liu2018robust}, etc. Owing to the increasing capability of computing power and the accelerated optimization algorithms, MPC has been utilized from the traditional process control to fields like complex dynamical systems \cite{kvasnica2019complexity,su2018self}, networked control systems \cite{li2017robust}, multi-agent systems \cite{li2016neighbor} and so on. It is worth to point out that for reducing the computational burden, the early termination is adopted in the optimization process to fulfill the real-time requirements, leading the solution to be inexact \cite{chen2018adaptive,kohler2019distributed}.

For large-scale systems, however, distributed MPC (DMPC) should be investigated to further reduce the computation resources compared with the centralized MPC. The most commonly used DMPC formulation is to describe the DMPC as a convex optimization problem that should be solved in a distributed fashion \cite{dunbar2012distributed,kohler2018distributed,liu2018distributed}, which can be referred to distributed consensus-based optimization \cite{shi2015extra,fazlyab2018distributed}. For given multi-agent systems, the goal of all the subsystems is to minimize a global objective function collaboratively without sharing the private objective function throughout the optimization process. Generally speaking, the distributed algorithms for convex optimization fall into three categories known as primal methods, dual methods and primal-dual methods. For primal methods, the convex optimization problem is solved in the primal domain with guaranteeing the consensus by introducing the penalty functions reflecting the disagreement \cite{patrascu2018convergence}. For dual methods, each subsystem solves the dual problem in a distributed manner to seek consensus \cite{yu2018convergence}. For the primal-dual methods, to obtain the saddle-point, the primal and dual problems are updated simultaneously associated with Lagrangian multipliers \cite{hale2017asynchronous}. It is worth to note that each of the methods has its advantages determining from the global objective function and the constraints formulations in the optimization problem. 

The coupled systems concerning DMPC scheme have received increasing attention recently as dynamical couplings are ubiquitous in practical applications. In \cite{dunbar2007distributed}, a distributed receding horizon control strategy was studied for a class of dynamically coupled nonlinear systems with consideration of decoupled constraints. In \cite{dai2017distributed}, the coupled probabilistic constraints were investigated in the context of distributed stochastic MPC. The compromise for satisfying the coupled constraints in a distributed way is to allow only one subsystem to optimize at each time instant, which is a widely utilized method. However, the distributed optimization problems have been rarely studied under the globally coupled constraints in the framework of DMPC. The limitation of the globally coupled constraints render the existing results on coupled systems cannot be generalized directly. The authors in \cite{wang2017distributed} presented a DMPC scheme for a group of discrete-time systems taking the local and global constraints into account. A dual problem was established to solve the DMPC optimization problem based on the Alternating Direction Multiplier Method (ADMM). It, however, is worth mentioning that the objective function considered in \cite{wang2017distributed} was with a quadratic form, which implies it cannot be directly extended to the scenario with a general nonlinear objective function. 

For nonlinear systems analysis, a well-known method is the contraction theory. By virtue of the fluid mechanics and differential geometry, the contraction theory has been first introduced in \cite{Lohmiller1998}. The traditional approach established by using Riemannian manifolds has been extended to many applications such as distributed nonlinear systems \cite{long2018distributed}, stochastic incremental systems \cite{pham2009contraction}, etc. Furthermore, some recent results were developed inspired by Finsler manifolds \cite{forni2014differential,chaffey2018control}. Nevertheless, it is worth noting that the existing results were obtained in the context of continuous-time dynamical systems. To the best of our knowledge, the contraction theory has been rarely investigated for discrete-time dynamical systems. Moreover, only a few results have been addressed using the contraction theory to analyze the convergence of the optimization algorithm \cite{nguyen2018contraction}. 

This paper formulates the DMPC optimization problem as a distributed consensus optimization problem (DCOP). To provide the capability of early termination with guaranteeing the convergence, a tightening constraint is constructed in the optimization. In addition, the primal-dual gradient optimization is adopted to solve the DCOP. The contraction theory based on Riemannian manifolds is utilized to analyze the convergence of the primal-dual gradient dynamics (PDGD). The contributions of this paper are mainly in three-fold as follows.

\begin{enumerate}[(1)]\label{contribution} 
	\item This paper investigates the DMPC subject to globally coupled constraints which can be formulated as a distributed consensus-based optimization problem. Furthermore, the inexact solver is taken into account to reduce unnecessary computations. The constraint tightening method is performed to allow premature termination with guaranteeing the convergence of the optimization process.
	
	\item The objective function in the DMPC optimization problem is considered as a nonlinear function in this paper. The dual problem is utilized to solve the DMPC optimization problem. Moreover, the local copies of the Lagrangian multiplier in the dual problem are introduced to achieve fully distributed. Thereafter, the primal-dual gradient dynamics are established to solve the consensus-based optimization problem. Owing to the tightening constraints, the local copies without being required to reach consensus but need to fulfill some specified bounds.
	
	\item Inspired by the contraction theory based on Riemannian manifolds, some sufficient conditions are given to guarantee the convergence of the PDGD. It is worth noting that the contraction theory is first used to analyze the optimization convergence in the context of discrete-time updating dynamics. In addition, the recursive feasibility and stability of the DMPC algorithm are rigorously analyzed.
	
\end{enumerate}

The remainder of this paper is organized as follows. Section II presents some necessary preliminaries adopted throughout this paper. In Section III, the considered optimization problem is formulated. The theoretical results are demonstrated in Section IV. Thereafter, Section V describes the proposed algorithm. We analyze the recursive feasibility and stability in Section VI. A numerical example is given in Section VII to verify the effectiveness of the proposed algorithm. Section VIII summarizes the paper.

The notations used in this paper are stated in the following. Denote real and natural number set as $\mathbb{R}$ and $\mathbb{N}$, respectively. The superscripts $\T$ and $+$ stand for transposition of a given matrix and the successor states. The subscripts of $\mathbb{N}_{\ge a}$ and $\mathbb{N}_{\left[b,\,c\right]}$ represent the integers in the intervals $\left[a,\,\infty\right)$ and $\left[b,c\right]$, respectively. Given two matrices $X$ and $Y$, the Kronecker product is denoted by $X\otimes Y$. Define the $P$-weighted norm of a given vector $x$ as $\left\| x\right\|_P = \sqrt{x^\T P x}$ with respect to the positive definite matrix $P$.

\section{Preliminaries}

\subsection{Graph Theory}

In this paper, the multi-agent system containing $l$ subsystems communicate with each other according to a weighted undirected graph denoted as $\mathcal{G} = \left(\mathcal{V},\mathcal{E},\mathcal{A}\right)$, where the set of agents is represented as $\mathcal{V}=\left\lbrace 1,2,\cdots,l\right\rbrace$, the set $\mathcal{E} = \mathcal{V}\times\mathcal{V}$ collects the undirected edges indicating the interconnections among subsystems, and the adjacent matrix $\mathcal{A}=[a_{ij}]_{l \times l}$ is symmetric with $a_{ii}=0$ implying that there is no self-edge in the graph, $a_{ij}=a_{ji}> 0$ denoting the undirected weight if $\left\lbrace i,j\right\rbrace\in\mathcal{E}$ and $a_{ij}= 0$ otherwise. Collect the neighbors of vertex $i$ in a set $\mathcal{N}_i = \left\lbrace j\in\mathcal{E}: a_{ij}> 0\right\rbrace$. The degree matrix is $\mathcal{D} =  \text{diag}\lbrace d_{i}\rbrace \in \mathbb{R}^{l\times l}$, where $d_{i} = \sum_{j\in\mathcal{N}_i}a_{ij}$. The Laplacian matrix is defined as $L = \mathcal{D}-\mathcal{A}$, and its eigenvalue decomposition is $L = Q\text{diag}\left(\nu_1,\nu_2,\cdots,\nu_M\right)Q^\T$, where $Q$ is an orthogonal matrix and $0=\nu_1<\nu_2\le\cdots\le\nu_l$. We can obtain $L = \sqrt{L}\sqrt{L}$ by introducing $\sqrt{L} = Q\text{diag}\left(\sqrt{\nu_1},\sqrt{\nu_2},\cdots,\sqrt{\nu_M}\right)Q^\T$. 

\subsection{Contraction Theory}

Some crucial definitions on Riemannian geometry are recapped for convergence analysis of the proposed algorithm by using contraction theory. For more details, please refer to \cite{Lohmiller1998,liu2017robust} and references therein. The Riemannian metric of two vectors $\delta_1, \delta_2$ on the tangent space of a given state manifold is a smoothly varying inner product $\left\langle \delta_1, \delta_2 \right\rangle_x = \delta_1^\T M\left(x\right) \delta_2$ with respect to a positive matrix function $M\left(x\right)$. Thus, we have $\|\delta\|_{x} = \sqrt{\left\langle \delta_1, \delta_2 \right\rangle_x}$. In this paper, we assume the matrix $M\left(x\right)=M$ to be constant. Given a pair of points $x\in\delta_1$ and $y\in \delta_2$, let the set of smooth curves connecting $x$ and $y$ be $\Gamma\left(x, y\right)$. There exists a piecewise smooth mapping $\gamma\in\Gamma\left(x, y\right)$ satisfying $\gamma\left(0\right) = x$ and $\gamma\left(1\right) = y$. Define the Riemannian length as $L\left(\gamma\right):=\int_{0}^{1} \|\gamma_s\|_{\gamma} ds$, the Riemannian energy $E\left(\gamma\right):=\int_{0}^{1} \|\gamma_s\|_{\gamma}^2 ds$, where $\gamma_s := \partial\gamma_s/\partial s$. Denote the Riemannian distance as $d\left(x,y\right) = \inf_{\gamma\in\Gamma\left(x,y\right)}L\left(\gamma\right)$. In this paper, we define $E\left(x,y\right) := d\left(x,y\right)^2$. 

Consider an autonomous nonlinear discrete-time system in the following form
\begin{equation}\label{autonomousSys}
\xi = f \left(\xi,k\right),
\end{equation} 
where $f$ is a smooth and differentiable function with respect to $\xi$. The differential dynamics can be denoted as $\delta \xi\left(k+1\right) = \frac{\partial f}{\partial x}\delta \xi\left(k\right)$. The optimal state trajectory is defined as a forward-complete solution of $\left(\ref{autonomousSys}\right)$. The Euler discretization of an exponentially controllable system possesses geometric convergence speed \cite{stuart1994numerical}. Thus, the optimal state trajectory $x^*$ is said to be global exponentially controllable if
\begin{equation}\label{locExpStable}
d\left(x\left(k\right),x^{*}\left(k\right)\right) \le C \tau ^{k}d\left(x\left(0\right),x^{*}\left(0\right)\right),\forall k>0,
\end{equation}
where the positive scalar $C$ and the convergence rate $\tau\in\left(0,1\right)$ are independent of the initial states. 

In accordance with \cite{Lohmiller1998}, the contraction region for the discrete-time system is similarly defined as follows. 

\textit{Definition 1:} Given a discrete-time system $x\left(k+1\right)=f\left(x\left(k\right),k\right)$, a region of state space is called a contraction region, if there exists a uniformly positive definite constant metric $M$, such that
\begin{equation}\label{contractionTheoryDef}
\frac{\partial f^{\T}}{\partial x}M\frac{\partial f}{\partial x} - M \le \left(\tau^2-1\right) M
\end{equation}
with the convergence rate $\tau\in\left(0,1\right)$.

\section{Problem Formulation}

Consider the multi-agent system containing $l$ subsystems under a weighted undirected graph. Each subsystem can be formulated as the following linear discrete-time dynamics:
\begin{equation}\label{sysForm}
x_i\left(t+1\right) = Ax_i\left(t\right)+Bu_i\left(t\right), \quad i=1,2,\cdots,l,
\end{equation}
with $x_i\left(t\right)\in\mathcal{X}_i$ and $u_i\left(t\right)\in\mathcal{U}_i$, where $x_i$ and $u_i$ are the state and input vectors of subsystem $i$, $\mathcal{X}_i\subset\mathbb{R}^n$ and $\mathcal{U}_i\subset\mathbb{R}^m$ are compact sets of the local constraints on state and input containing the origin as their inner point, respectively. Moreover, the following globally coupled constraint is taken into account:
\begin{equation}\label{coupledConstraint}
\sum_{i=1}^{l}\Phi_i^x x_i\left(t\right)+\Phi_i^u u_i\left(t\right)\le \boldsymbol{1}_p, \quad \forall t\ge 0,
\end{equation}
where $\Phi_i^x\in\mathbb{R}^{p\times n}$, $\Phi_i^u\in\mathbb{R}^{p\times m}$ are the matrices defining the coupled constraint, and $\boldsymbol{1}_p$ represents the $p$-vector with all ones. 

%The objective function for each subsystem is designed as
%\begin{equation}\label{objFcn}
%\begin{aligned}
%J_i\left(x_i\left(\cdot|t\right),u_i\left(\cdot|t\right)\right)&:=\sum_{s=0}^{N-1}\left(\left\|x_i\left(s|t\right) \right\|_{Q_i}^2+\left\|u_i\left(s|t\right) \right\|_{R_i}^2\right)\\
%&\qquad\qquad +\left\|x_i\left(N|t\right) \right\|_{P_i}^2
%\end{aligned}
%\end{equation}
%where $N\in\mathbb{N}_{>0}$ is the prediction horizon, $Q_i> 0$, $R_i> 0$, $P_i> 0$ stand for the weight matrices, $x_i\left(s|t\right)$ and $u_i\left(s|t\right)$ depict the state and input sequences for $s\in\mathbb{N}_{\left[1,N-1\right]}$, respectively. Thus, the standard MPC optimization problem is formulated by
For subsystem $i$ at time instant $t$, define the objective function $J_i\left(x_i,u_i\right):\mathbb{R}^n\times \mathbb{R}^m\rightarrow \mathbb{R}$ as  
\begin{equation}\label{ObjFcn}
J_i\left(x_i,u_i\right) = \sum_{s=0}^{N} F\left(x_i\left(s|t\right),u_i\left(s|t\right)\right)+V_f\left(x_i\left(N|t\right)\right)
\end{equation}
where $N$ is the prediction horizon. The following assumptions are crucial for theoretical analysis in this paper.

\textit{Assumption 1:} \cite{fazlyab2018distributed} The local objective function $J_i\left(x_i,u_i\right)$ is twice-continuously differentiable and $m_J$-strongly convex with Lipschitz gradient $L_J$ with respect to $u_i$, ie., for all $u_i,u_i'\in\mathcal{U}_i$,
\begin{subequations}\label{assumption1}
\begin{align}
J_i\left(x_i,u_i\right) &\ge J_i\left(x_i,u_i'\right)+ \frac{m_J}{2}\left\|u_i-u_i' \right\|_2^2 \nonumber\\
&+ \frac{m_J}{2}\left(\nabla_{u_i} J_i\left(x_i,u_i'\right)\right)^\T\left(u_i'-u_i\right),\\
J_i\left(x_i,u_i\right) &\le J_i\left(x_i,u_i'\right)+ \frac{L_J}{2}\left\|u_i-u_i' \right\|_2^2 \nonumber\\
&+ \left(\nabla_{u_i} J_i\left(x_i,u_i'\right)\right)^\T\left(u_i-u_i'\right).
\end{align}
\end{subequations} 

\textit{Remark 1:} Assumption 1 is standard in convergence analysis of convex optimization. We can rewrite Eq. $\left(\ref{assumption1} \right)$ into another form as
\begin{equation}\label{assumption1_trans}
m_J \boldsymbol{I} \le \nabla_{u_i}^2 J_i\left(x_i,u_i\right) \le L_J \boldsymbol{I}
\end{equation}
with $0<m_J\le L_J<\infty$.

\textit{Assumption 2:} \cite{qu2018exponential} Under Assumption 1, for any $\boldsymbol{u}_i$ and the optimal $u^*$, there exists an invertible symmetric matrix $H_i\left(u_i\right)$ satisfying $m_J\boldsymbol{I}\le H_i\left(u_i\right)\le L_J\boldsymbol{I}$ such that
\begin{equation}\label{assumption2}
\nabla J_i\left(x_i,u_i\right)-\nabla J_i\left(x_i,u_i^*\right) = H_i\left(u_i\right) \left(u_i-u_i^* \right).
\end{equation}

\textit{Remark 2:} Notice that $H_i\left(u_i\right)$ is a time-varying matrix with respect to $u_i$. For convenience, we use $H_i$ for short in this paper.

Inspired by \cite{rubagotti2014stabilizing,wang2018accelerated}, we define
\begin{equation}\label{U_x}
\begin{aligned}
\tilde{\mathcal{U}}_i\left(x_i\left(t\right)\right) := &\left\lbrace \tilde{\boldsymbol{u}}_i: \tilde{\boldsymbol{x}}_i\left(s+1|t\right) = A \tilde{\boldsymbol{x}}_i\left(s|t\right)\right.\\
&\left.+B \tilde{\boldsymbol{u}}_i\left(s|t\right), \tilde{\boldsymbol{x}}_i\left(0|t\right)=x_i\left(t\right),\tilde{\boldsymbol{x}}_i\left(s|t\right)\in\mathcal{X}_i,\right.\\
&\left.\tilde{\boldsymbol{u}}_i\left(s|t\right)\in\mathcal{U}_i, \tilde{\boldsymbol{x}}_i\left(N|t\right)\in\Omega_i^f,s\in\mathbb{Z}_{\left[0,N-1\right]}\right\rbrace,
\end{aligned}
\end{equation} 
where $\tilde{\boldsymbol{x}}_i\left(s|t\right)$ and $\tilde{\boldsymbol{u}}_i\left(s|t\right)$ depict the predicted state and input sequences for $s\in\mathbb{N}_{\left[0,N-1\right]}$, respectively, $N\in\mathbb{Z}_{>0}$ represents the prediction horizon and $\Omega_i^f := \left\lbrace x_i|V_f\left(x_i\right)\le\eta_f \right\rbrace$ denotes the terminal set. In addition, the following assumption is given for stability analysis.

\textit{Assumption 3:} \cite{hashimoto2016self} Given a positive scalar $\eta>\eta_f$, there exists a local state-feedback control gain $K$ such that $u_i = Kx_i\in\mathcal{U}_i$. Moreover, it holds that
\begin{equation}\label{assumption3}
V_f\left(A_Kx_i\left(t\right)\right)-V_f\left(x_i\left(t\right)\right)\le-F\left(x_i\left(t\right),Kx_i\left(t\right)\right)
\end{equation}
with $A_K = A+BK$ for all $x_i\in\Omega_i:=\left\lbrace x_i|V_f\left(x_i\right)\le\eta \right\rbrace$.

It is worth to clarify that the terminal set used in $\left(\ref{U_x} \right)$ is smaller than which is defined in Assumption 3. In addition, the following condition should be satisfied for all $x_i\left(t\right)\in\Omega_i^f$ to fulfill the globally coupled constraint in $\left(\ref{coupledConstraint} \right)$.
\begin{equation}\label{localCoupleConstraint}
\sum_{i=1}^{l}\left(\Phi_i^x +\Phi_i^u K_i\right)x_i\left(s|t\right)\le \boldsymbol{1}_p, \quad \forall s\in\mathbb{Z}_{>0}.
\end{equation}

In what follows, we formulate the standard MPC optimization problem as
\begin{subequations}\label{standardMPC}
	\begin{align}	
	u_i^*\left(s|t\right) = \arg \min_{\tilde{\boldsymbol{u}}_i\left(s|t\right)} \sum_{i=1}^{l}J_i\left(\tilde{\boldsymbol{x}}_i\left(s|t\right),\tilde{\boldsymbol{u}}_i\left(s|t\right)\right)
	\end{align}
	\begin{align}
	s.t.\quad&\tilde{\boldsymbol{u}}_i\left(s|t\right)\in\tilde{\mathcal{U}}_i\left(x_i\left(t\right)\right),\\
	&\sum_{i=1}^{l}\Phi_i^x \tilde{\boldsymbol{x}}_i\left(s|t\right)+\Phi_i^u \tilde{\boldsymbol{u}}_i\left(s|t\right)\le \boldsymbol{1}_p,\label{standardMPC:coupled}
	\end{align}
\end{subequations}
with $ s\in\mathbb{Z}_{\left[0,N-1\right]}$, where $\tilde{\mathcal{U}}_i\left(x_i\left(t\right)\right)$ is defined in $\left(\ref{U_x}\right)$ and $\left(\ref{standardMPC:coupled}\right)$ implying the satisfaction of the globally coupled constraint in $\left(\ref{coupledConstraint}\right)$.

To accelerate the optimization, one efficient way is to allow premature termination of the optimization problem in $\left(\ref{standardMPC} \right)$ with guaranteeing the convergence, which can avoid unnecessary computations. It is worth mentioning that the early termination may result in errors even infeasibility of the optimization problem. Taking the inexactness of optimization solver into consideration, the $\epsilon$-strict feasibility is investigated by introducing the tightening constraint in this paper.

\textit{Definition 2:} \cite{kohler2019distributed} Given a polytopic constraint as $f\left(x\right)\le b$, the $\epsilon$-strictly feasible solution is the vector $x$ which satisfies $f\left(x\right)\le b+\epsilon \boldsymbol{1}$ with $\boldsymbol{1}$ in proper dimension.

 Inspired by the tightening constraint method, we can transform the global coupled constraint in $\left(\ref{standardMPC:coupled} \right)$ as follows
\begin{equation}\label{tighteningConstraint}
\sum_{i=1}^{l}\Phi_i^x \tilde{\boldsymbol{x}}_i\left(s|t\right)+\Phi_i^u \tilde{\boldsymbol{u}}_i\left(s|t\right)\le \left(1-\epsilon l \left(s+1\right)\right)\boldsymbol{1}_p,
\end{equation}
with $ s\in\mathbb{Z}_{\left[0,N-1\right]}$ where $\epsilon\in\mathbb{R}_{\left(0,\frac{1}{Nl}\right)}$ is a user-defined parameter reflecting the tolerance to the violation of the coupled constraint. For convenience, we rewrite $\left(\ref{tighteningConstraint} \right)$ as a standard form in optimization. 
\begin{equation}\label{tighteningConstraintOpt}
\sum_{i=1}^{l} f_i\left(x_i\left(t\right),\tilde{\boldsymbol{u}}_i\left(s|t\right)\right)\le b\left(\epsilon \right),
\end{equation}
with $ s\in\mathbb{Z}_{\left[0,N-1\right]}$, where 
\begin{equation}\label{tighteningConstraintOpt1}
f_i\left(x_i\left(s|t\right),\tilde{\boldsymbol{u}}_i\left(s|t\right)\right):=F_i x_i\left(s|t\right)+ G_i \tilde{\boldsymbol{u}}_i\left(s|t\right)
\end{equation} 
with $F_i\in\mathbb{R}^{Np\times n}$ and $G_i\in\mathbb{R}^{Np\times Nm}$ being appropriate matrices obtained from $\Phi_i^x$ and $\Phi_i^u$ by expressing $\tilde{\boldsymbol{x}}_i\left(s|t\right)$ in terms of $x_i\left(t\right)$ and $\tilde{\boldsymbol{u}}_i\left(s|t\right)$, and 
\begin{equation}\label{tighteningConstraintOpt2}
b\left(\epsilon \right) := \left[\left(1-\epsilon l\right)\boldsymbol{1}_p^\T, \left(1-2\epsilon l\right)\boldsymbol{1}_p^\T,\cdots,\left(1-\epsilon Nl\right)\boldsymbol{1}_p^\T \right]^\T.
\end{equation} 

The following assumption is established based on the linear independence property of the coupled constraint.

\textit{Assumption 4:} There exist two scalars $\underline{\zeta}>0$ and $\overline{\zeta}>0$ such that the full row rank matrix $F_i$ satisfying
\begin{equation}\label{assumption4}
\underline{\zeta}\boldsymbol{I} \le G_i^\T G_i \le \overline{\zeta}\boldsymbol{I}
\end{equation}
where $\boldsymbol{I}$ is with appropriate dimension.

Moreover, the coupled constraint in the terminal set can be rewritten as the following form, correspondingly.
\begin{equation}\label{terminalCOupledConstraint}
\sum_{i=1}^{l}\left(\Phi_i^x +\Phi_i^u K_i\right)x_i\left(s|t\right)\le \left(1-\epsilon Nl\right)\boldsymbol{1}_p,
\end{equation}
with $s\in\mathbb{Z}_{>0}$ for $x_i\in\Omega_i^f$. By utilizing the constraint in $\left(\ref{tighteningConstraintOpt} \right)$, the DMPC optimization problem in $\left(\ref{standardMPC} \right)$ can be rewritten as follows
\begin{subequations}\label{erMPC}
	\begin{align}	
	\tilde{\boldsymbol{u}}_i^*\left(s|t\right) = \arg \min_{\tilde{\boldsymbol{u}}_i\left(t\right)} \sum_{i=1}^{l}J_i\left(x_i\left(t\right),\tilde{\boldsymbol{u}}_i\left(s|t\right)\right)\label{erMPC:ObjFcn}
	\end{align}
	\begin{align}
s.t.\quad	&\tilde{\boldsymbol{u}}_i\left(s|t\right)\in\tilde{\mathcal{U}}_i\left(x_i\left(t\right)\right),\label{erMPC:U}\\
	&\sum_{i=1}^{l} f_i\left(x_i\left(t\right),\tilde{\boldsymbol{u}}_i\left(s|t\right)\right)\le b\left(\epsilon \right),\label{erMPC:coupled}
	\end{align}
\end{subequations}
with $s\in\mathbb{Z}_{\left[0,N-1\right]}$, where $\tilde{\mathcal{U}}_i\left(x_i\right)$ is defined in $\left(\ref{U_x} \right)$.

\section{Main Results}

\subsection{The Dual Form}

To solve the optimization problem in $\left(\ref{erMPC} \right)$, the dual form is considered in the following. The Lagrangian of optimization problem $\left(\ref{erMPC}\right)$ is given as
\begin{equation}\label{lagrangian1}
\mathcal{L}\left(x_i,\tilde{\boldsymbol{u}}_i,\boldsymbol{\lambda}\right) = \sum_{i=1}^{l}J_i\left(x_i,\tilde{\boldsymbol{u}}_i\right)+\boldsymbol{\lambda}^\T\left(\sum_{i=1}^{l} f_i\left(x_i,\tilde{\boldsymbol{u}}_i\right)- b\left(\epsilon \right)\right)
\end{equation}
for $\tilde{\boldsymbol{u}}_i\in\tilde{\mathcal{U}}_i\left(x_i\right)$, where $\boldsymbol{\lambda}\in\mathbb{R}^{Np}$ is the Lagrangian multiplier. Thus, the dual problem can be described as
\begin{equation}\label{dualProblem1}
\max_{\boldsymbol{\lambda}\ge 0}\;\min_{\tilde{\boldsymbol{u}}_i\in\tilde{\mathcal{U}}_i\left(x_i\right)}  \;\mathcal{L}\left(x_i,\tilde{\boldsymbol{u}}_i,\boldsymbol{\lambda}\right)
\end{equation}
which is equivalent to
\begin{equation}\label{dualProblem2}
\min_{\boldsymbol{\lambda}\ge 0}\;\max_{\tilde{\boldsymbol{u}}_i\in\tilde{\mathcal{U}}_i\left(x_i\right)}  \; -\mathcal{L}\left(x_i,\tilde{\boldsymbol{u}}_i,\boldsymbol{\lambda}\right).
\end{equation}
The Karush-Kuhn-Tucker (KKT) conditions for depicting the optimal pair $\left(\tilde{\boldsymbol{u}}_i^*,\boldsymbol{\lambda}^* \right)$ are
\begin{subequations}\label{dualKKT}
\begin{align}
\nabla_{\tilde{\boldsymbol{u}}_i^*} \mathcal{L}\left(\tilde{\boldsymbol{u}}_i^*,\boldsymbol{\lambda}^*\right)&= \nabla_{\tilde{\boldsymbol{u}}_i^*}\sum_{i=1}^{l}J_i\left(x_i,\tilde{\boldsymbol{u}}_i^*\right)+G_i^\T\boldsymbol{\lambda}^*=0,\label{dualKKT:1}\\
\nabla_{\boldsymbol{\lambda}^*} \mathcal{L}\left(\tilde{\boldsymbol{u}}_i^*,\boldsymbol{\lambda}^*\right) &= \sum_{i=1}^{l} f_i\left(x_i,\tilde{\boldsymbol{u}}_i^*\right)- b\left(\epsilon \right)=0.\label{dualKKT:2}
\end{align}
\end{subequations}

\textit{Remark 3:} For each subsystem, we rewrite the dual problem in $\left(\ref{dualProblem2} \right)$ in the following form
\begin{equation}\label{dualProblem3}
\min_{\boldsymbol{\lambda}\ge 0} \; \Psi_i\left(\boldsymbol{\lambda}\right)
\end{equation}
where 
\begin{equation}\label{dualProblem4}
\Psi_i\left(\boldsymbol{\lambda}\right):=\max_{\tilde{\boldsymbol{u}}_i\in\tilde{\mathcal{U}}_i\left(x_i\right)} -J_i\left(x_i,\tilde{\boldsymbol{u}}_i\right)-\boldsymbol{\lambda}^\T\left(f_i\left(x_i,\tilde{\boldsymbol{u}}_i\right)- \frac{b\left(\epsilon \right)}{l}\right).
\end{equation}
By Dadnskin's theorem \cite{bertsekas1999nonlinear}, the dual gradient is obtained as $\nabla\Psi_i\left(\boldsymbol{\lambda}\right) = -\left(f_i\left(x_i,\tilde{\boldsymbol{u}}_i\left(\boldsymbol{\lambda}\right)\right)- \frac{b\left(\epsilon \right)}{l}\right)$.

\subsection{Distributed Optimization with Laplacian Consensus}

It is worth to point out that the Lagrangian multiplier $\boldsymbol{\lambda}$ in $\left(\ref{dualProblem3} \right)$ is a global variable such that the optimization problem cannot be solved in a distributed way. Resorting to the Laplacian consensus, the optimization problem in $\left(\ref{dualProblem3} \right)$ can be transformed into
\begin{subequations}\label{consensusOpt}
	\begin{align}	
	\min_{\boldsymbol{\lambda}_i\ge 0} \Psi_i\left(\boldsymbol{\lambda}_i\right)
	\end{align}
	\begin{align}
	s.t.\quad\sqrt{\boldsymbol{L}}\boldsymbol{\Lambda} = 0,\label{consensusOpt:constraint}
	\end{align}
\end{subequations}
where $\sqrt{\boldsymbol{L}}:=\sqrt{L}\otimes \boldsymbol{I}_{Np}$ with $L$ being the Laplacian matrix corresponding to the topology of the communication graph, $\boldsymbol{\lambda}_i$ is the local copy of $\boldsymbol{\lambda}$ for the subsystem $i$, and $\boldsymbol{\Lambda} := \left[\boldsymbol{\lambda}_1^\T,\boldsymbol{\lambda}_2^\T,\cdots, \boldsymbol{\lambda}_l^\T \right]^\T$ is the vector which stacks the local Lagrangian multipliers $\boldsymbol{\lambda}_i$. For the individual subsystem, we can obtain the Lagrangian of optimization problem in $\left(\ref{consensusOpt} \right)$ as follows
\begin{equation}\label{lagrangian2}
\mathcal{D}_i\left(\boldsymbol{\lambda}_i,\boldsymbol{\mu}_i\right) = \Psi_i\left(\boldsymbol{\lambda}_i\right)+\boldsymbol{\mu}_i^\T \sqrt{\boldsymbol{L}}\boldsymbol{\Lambda} ,
\end{equation}
where $\boldsymbol{\mu}_i\in\mathbb{R}^{Np}$ is the Lagrangian multiplier of $i$th subsystem. The KKT conditions for the optimal pair $\left(\boldsymbol{\lambda}_i^*,\boldsymbol{\mu}_i^*\right)$ are described as
\begin{subequations}\label{KKT}
\begin{align}
	\nabla_{\boldsymbol{\lambda}_i^*} \mathcal{D}_i\left(\boldsymbol{\lambda}_i^*,\boldsymbol{\mu}_i^*\right)&= \nabla\Psi_i\left(\boldsymbol{\lambda}_i^*\right)+ \sqrt{\boldsymbol{L}}\boldsymbol{\mu}_i^*=0,\\
	\nabla_{\boldsymbol{\mu}_i^*} \mathcal{D}_i\left(\boldsymbol{\lambda}_i^*,\boldsymbol{\mu}_i^*\right) &= \sqrt{\boldsymbol{L}}\boldsymbol{\Lambda}^*=0.
\end{align}
\end{subequations}

\subsection{Primal-Dual Gradient Dynamics}

The primal-dual gradient method is adopted to solve the optimization problems in $\left(\ref{dualProblem3} \right)$ and $\left(\ref{consensusOpt} \right)$. The primal-dual gradient dynamics are given as
\begin{subequations}\label{PDGD1}
	\begin{align}
	\boldsymbol{\lambda}_i^{k+1}-\boldsymbol{\lambda}_i^{k}&= - \alpha_i\nabla_{\boldsymbol{\lambda}_i^k} \mathcal{D}_i\left(\boldsymbol{\lambda}_i^k,\boldsymbol{\mu}_i^k\right)\nonumber\\
	&= - \alpha_i\nabla\Psi_i\left(\boldsymbol{\lambda}_i^k\right)- \alpha_i\sqrt{\boldsymbol{L}}\boldsymbol{\mu}_i^k, \\
	\boldsymbol{\mu}_i^{k+1}-\boldsymbol{\mu}_i^k &= \beta_i\nabla_{\boldsymbol{\mu}_i^k} \mathcal{D}_i\left(\boldsymbol{\lambda}_i^k,\boldsymbol{\mu}_i^k\right) = \beta_i\sqrt{\boldsymbol{L}}\boldsymbol{\Lambda}^k,
	\end{align}
\end{subequations}
where $\alpha>0$ and $\beta>0$ are the step-sizes, $\nabla \Psi_i\left(\boldsymbol{\lambda}_i^k\right) = -\left(f_i\left(x_i,\tilde{\boldsymbol{u}}_i\left(\boldsymbol{\lambda}_i^k\right)\right)- \frac{b\left(\epsilon \right)}{l}\right)$ with
\begin{equation}\label{u_i}
\tilde{\boldsymbol{u}}_i^k = \arg\min_{\tilde{\boldsymbol{u}}_i}J_i\left(\tilde{\boldsymbol{u}}_i\right)+\left(\boldsymbol{\lambda}_i^k\right)^\T\left(f_i\left(x_i,\tilde{\boldsymbol{u}}_i\right)-\frac{b\left(\epsilon \right)}{l}\right).
\end{equation}

\textit{Remark 4:} It is worth to note that the Laplacian matrix $L$ cannot be locally computed such that Eq. $\left(\ref{PDGD1}\right)$ is not distributed. We introduce a new variable to scale the Lagrangian multiplier $\boldsymbol{\mu}_i$ as $\boldsymbol{\gamma}_i =\sqrt{\boldsymbol{L}} \boldsymbol{\mu}_i$. Thus, the PDGD in $\left(\ref{PDGD1} \right)$ is equivalent to
\begin{subequations}\label{PDGD2}
	\begin{align}
	\boldsymbol{\lambda}_i^{k+1}-\boldsymbol{\lambda}_i^k&= -\alpha_i\nabla\Psi_i\left(\boldsymbol{\lambda}_i^k\right)-\alpha_i\boldsymbol{\gamma}_i^k\label{PDGD2:lambda}, \\
	\boldsymbol{\gamma}_i^{k+1}-\boldsymbol{\gamma}_i^k &=  \beta_i \boldsymbol{L} \boldsymbol{\Lambda}^k=  \beta_i d_i\boldsymbol{\lambda}_i^k-\beta_i\sum_{j\in\mathcal{N}_i}a_{ij}\boldsymbol{\lambda}_j^k\label{PDGD2:gamma}.
	\end{align}
\end{subequations} 
Notice that from a global perspective, the introduced $\boldsymbol{\gamma}^0 = \left[\boldsymbol{\gamma}_{1}^0,\boldsymbol{\gamma}_{2}^0,\cdots,\boldsymbol{\gamma}_{l}^0 \right]^\T = \sqrt{\boldsymbol{L}}\boldsymbol{\Lambda}^0$ requires the initial values to satisfy $\sum_{i=1}^{l}\boldsymbol{\gamma}_{i}^0 = 0$.

%The PDGD iteration for the $i$th subsystem can be described in the following
%\begin{subequations}\label{PDGD3}
%	\begin{align}
%	\lambda_i^+ -\lambda_i&= -\nabla_{\lambda_i}\Psi_i\left(\lambda_i\right)-\gamma_i \\
%	\gamma_i^+ - \gamma_i &= \tau \sum_{j\in\mathcal{N}_i}a_{ij}\left(\lambda_i-\lambda_j\right)
%	\end{align}
%\end{subequations} 

\subsection{Contraction of the PDGD}

Inspired by the Riemannian geometry, some sufficient conditions guaranteeing the exponential convergence of the PDGD in $\left(\ref{PDGD2} \right)$ are given in the following theorem. 

\textit{Theorem 1:} Under Assumptions 2 and 4, the PDGD in $\left(\ref{PDGD2} \right)$ has the exponential convergence rate $\tau =  \sqrt{1-\rho\alpha_i\beta_i d_i}$ with Riemannian metric
\begin{equation}\label{RieMatrix}
M = 
\begin{bmatrix}
\beta_i d_i\boldsymbol{I}& \alpha_i\beta_i d_i \boldsymbol{I}\\
\alpha_i\beta_i d_i\boldsymbol{I} & \alpha_i \boldsymbol{I}
\end{bmatrix},
\end{equation}
if 
\begin{subequations}\label{alpha_beta}
\begin{align}
\alpha_i&\le\frac{\left(1-\rho\right)\underline{\sigma}}{2\overline{\sigma}^2-\left(3+\rho \right)\underline{\sigma}\beta_i d_i},\\
\beta_i&\le\frac{\overline{\sigma}}{2 d_i},
\end{align}
\end{subequations}
where $\underline{\sigma}$ and $\overline{\sigma}$ derived from $\left(\ref{assumption2} \right)$ and $\left(\ref{assumption4} \right)$ are two proper scalars satisfying $0<\underline{\sigma}\le  \overline{\sigma} <\infty$, $0<\rho<\min\left(1, \frac{1}{\alpha_i\beta_i d_i}\right)$ is a user-defined scalar parameter and $\boldsymbol{I}$ is with appropriate dimension.

\textit{Proof.} Stack the Lagrangian multipliers $\boldsymbol{\lambda}_i$ and $\boldsymbol{\gamma}_i$ into a vector $y = \left[\boldsymbol{\lambda}_i^\T, \boldsymbol{\gamma}_i^\T \right]^\T$. The optimal vector can be similarly defined as $y^* = \left[\left(\boldsymbol{\lambda}_i^*\right)^\T, \left(\boldsymbol{\gamma}_i^*\right)^\T \right]^\T$. According to the updating dynamics of the primal-dual gradient in $\left(\ref{PDGD2} \right)$, we can obtain the differential dynamics of PDGD as follows
\begin{equation}\label{diffDynPDGD}
\begin{aligned}
&y^+-\left(y^*\right)^+\\
=&\begin{bmatrix}
\boldsymbol{\lambda}_i-\alpha_i\nabla\Psi_i\left(\boldsymbol{\lambda}_i\right)-\alpha_i\boldsymbol{\gamma}_i-\left(\boldsymbol{\lambda}_i^*-\alpha_i\nabla\Psi_i\left(\boldsymbol{\lambda}_i^*\right)-\alpha_i\boldsymbol{\gamma}_i^*\right)\\
\boldsymbol{\gamma}_i+\beta_i d_i\boldsymbol{\lambda}_i -\left(\boldsymbol{\gamma}_i^*+\beta_i d_i\boldsymbol{\lambda}_i^*\right) 
\end{bmatrix}\\
=&\begin{bmatrix}
\left(1-\alpha_i \Theta_i\right)\left(\boldsymbol{\lambda}_i-\boldsymbol{\lambda}_i^*\right)-\alpha_i\left(\boldsymbol{\gamma}_i-\boldsymbol{\gamma}_i^*\right)\\
\left(\boldsymbol{\gamma}_i-\boldsymbol{\gamma}_i^*\right)+\beta_i d_i\left(\boldsymbol{\lambda}_i -\boldsymbol{\lambda}_i^*\right) 
\end{bmatrix}= \Xi \left(y-y^*\right),
\end{aligned}
\end{equation}
where $\Theta_i = G_i H_i^{-1}G_i^\T$ and 
\begin{equation}\label{Xi}
\Xi = 
\begin{bmatrix}
\boldsymbol{I}- \alpha_i \Theta_i & -\alpha_i\boldsymbol{I}\\
 \beta_i d_i\boldsymbol{I}& \boldsymbol{I}
\end{bmatrix}.
\end{equation}
The second equality of $\left(\ref{diffDynPDGD} \right)$ follows from the condition in $\left(\ref{u_i} \right)$ and the KKT condition in $\left(\ref{dualKKT:1} \right)$ under Assumption 2. By constructing the Riemannian metric as $\left(\ref{RieMatrix} \right)$, the difference of the Riemannian energy between the adjacent updates can be written as
\begin{equation}\label{diffRieDist}
\begin{aligned}
&E\left(y^+,\left(y^*\right)^+\right)-E\left(y,y^*\right)\\
=&\left\| y^+-\left(y^*\right)^+\right\|_{M}^2-\left\| y-y^*\right\|_{M}^2\\
=&\left(y-y^*\right)^\T \left(\Xi^\T M\Xi - M  \right) \left(y-y^*\right)\\
=&\left(y-y^*\right)^\T\Pi\left(y-y^*\right),
\end{aligned}
\end{equation}
where 
\begin{equation}\label{Pi}
\Pi = 
\begin{bmatrix}
\Pi_{1} & \Pi_{2}^\T\\
\Pi_{2} & \Pi_{3}
\end{bmatrix}
\end{equation}
with
\begin{subequations}
	\begin{align}
	\Pi_{1} &= \alpha_i\beta_i^2 d_i^2\left(\left(\boldsymbol{I}-\alpha_i \Theta_i\right)^\T+\left(\boldsymbol{I}-\alpha_i \Theta_i\right)+\boldsymbol{I}\right)\nonumber\\
	&\quad +\beta_i d_i \left(\left(\boldsymbol{I}-\alpha_i \Theta_i\right)^\T\left(\boldsymbol{I}-\alpha_i \Theta_i\right)-\boldsymbol{I}\right),\\
	\Pi_{2} &=  -\alpha_i^2\beta_i^2 d_i^2 \boldsymbol{I},\\
	\Pi_{3} &= -\alpha_i^2\beta_i d_i \boldsymbol{I}.
	\end{align}
\end{subequations}
Therefore, to prove $\left(\ref{contractionTheoryDef} \right)$, it suffices to prove that $\Pi\le\left(\tau^2-1\right) M$. Letting $\Upsilon = \left(\tau^2-1\right) M-\Pi$, we can obtain
\begin{equation}\label{Upsilon}
\Upsilon = 
\begin{bmatrix}
\Upsilon_{1}&\Upsilon_{2}^\T\\
\Upsilon_{2}&\Upsilon_{3}
\end{bmatrix}
\end{equation}
with
\begin{subequations}
	\begin{align}
	\Upsilon_{1} &= \tau^2\beta_i d_i\boldsymbol{I} -\beta_i d_i \left(\boldsymbol{I}-\alpha_i \Theta_i\right)^\T\left(\boldsymbol{I}-\alpha_i \Theta_i\right)\nonumber\\
	& \quad-\alpha_i\beta_i^2 d_i^2\left(\left(\boldsymbol{I}-\alpha_i\Theta_i\right)^\T+\left(\boldsymbol{I}-\alpha_i \Theta_i\right)+\boldsymbol{I}\right),\\
	\Upsilon_{2} &= \left(\tau^2 \alpha_i\beta_i d_i - \alpha_i\beta_i d_i + \alpha_i^2\beta_i^2 d_i^2\right) \boldsymbol{I},\\
	\Upsilon_{3} &= \left(\tau^2 \alpha_i -\alpha_i+ \alpha_i^2\beta_i d_i\right) \boldsymbol{I}.
	\end{align}
\end{subequations}
Resorting to the Schur's complement, to prove $\Upsilon\ge0$, it is sufficient to prove $\Upsilon_{3}>0$ and $\Upsilon_{1}-\Upsilon_{2}\Upsilon_{3}^{-1}\Upsilon_{2}^\T\ge0$. By introducing a user-defined parameter $0<\rho<\min\left(1, \frac{1}{\alpha_i\beta_i d_i}\right)$ such that $\tau = \sqrt{1-\rho\alpha_i\beta_i d_i}$, one can get
\begin{equation}\label{Upsilon3}
\Upsilon_{3} = \left(1-\rho \right)\alpha_i^2\beta_i d_i\boldsymbol{I} > 0
\end{equation}
By using $\alpha_i<\frac{\left(1-\rho\right)\underline{\sigma}}{2\overline{\sigma}^2-\left(3+\rho \right)\underline{\sigma}\beta_i d_i}$ and $\beta_i<\frac{\overline{\sigma}}{2d_i}$, one have
\begin{equation}\label{Upsilon1}
\begin{aligned}
\Upsilon_{1}-&\Upsilon_{2}^\T\Upsilon_{3}^{-1}\Upsilon_{2}\\
=& \left(1-\rho\alpha_i\beta_i d_i\right)\beta_i d_i\boldsymbol{I} -\beta_i d_i \left(\boldsymbol{I}-\alpha_i\Theta_i\right)^\T\left(\boldsymbol{I}-\alpha_i \Theta_i\right) \\
&-\alpha_i\beta_i^2 d_i^2\left(\left(\boldsymbol{I}-\alpha_i \Theta_i\right)^\T+\left(\boldsymbol{I}-\alpha_i \Theta_i\right)+\boldsymbol{I}\right)\\
&-\frac{\left(\left(1-\rho\alpha_i\beta_i d_i\right) \alpha_i\beta_i d_i - \alpha_i\beta_i d_i + \alpha_i^2\beta_i^2 d_i^2\right)^2}{\left(\left(1-\rho\alpha_i\beta_i d_i\right) \alpha_i -\alpha_i+ \alpha_i^2\beta_i d_i\right)}\boldsymbol{I}\\
=&-\left(3+\rho \right)\alpha_i\beta_i^2 d_i^2\boldsymbol{I}-\left(1-\rho \right)\alpha_i^2\beta_i^3 d_i^3\boldsymbol{I}-\alpha_i^2\beta_i d_i \Theta_i^\T \Theta_i\\
&+\left(\alpha_i^2\beta_i^2 d_i^2+\alpha_i\beta_i d_i\right)\left(\Theta_i^\T+ \Theta_i\right)\\
\ge&\left(-\left(3+\rho \right)\alpha_i\beta_i^2 d_i^2-\left(1-\rho \right)\alpha_i^2\beta_i^3 d_i^3-\alpha_i^2\beta_i d_i \overline{\sigma}^2\right.\\
&\left.+2\underline{\sigma}\left(\alpha_i^2\beta_i^2 d_i^2+\alpha_i\beta_i d_i\right)\right)\boldsymbol{I}\ge 0
\end{aligned}
\end{equation}
In light of Assumptions $2$ and $4$, it suffices to have
\begin{equation}\label{sigmaRange}
\underline{\sigma}\boldsymbol{I}  \le \Theta_i \le \overline{\sigma}\boldsymbol{I} 
\end{equation}
where $\underline{\sigma}$ and $\overline{\sigma}$ are two scalars satisfying $0<\underline{\sigma}\le  \overline{\sigma} <\infty$, which are derived from $\left(\ref{assumption2} \right)$ and $\left(\ref{assumption4} \right)$ appropriately. Thus, we can obtain the last inequality of $\left(\ref{Upsilon1} \right)$ following from
\begin{subequations}\label{Upsilon1:1}
\begin{align}
&\frac{1}{2}\underline{\sigma}\left(\left(3+\rho \right)\alpha_i\beta_i d_i+\left(1-\rho \right)\alpha_i^2\beta_i^2 d_i^2\right)\boldsymbol{I}\nonumber\\
&\qquad\ge\left(\left(3+\rho \right)\alpha_i\beta_i^2 d_i^2+\left(1-\rho \right)\alpha_i^2\beta_i^3 d_i^3\right)\boldsymbol{I}\\
&\frac{1}{2}\underline{\sigma}\left(\left(3+\rho \right)\alpha_i^2\beta_i^2 d_i^2+\left(1-\rho\right)\alpha_i\beta_i d_i\right)\ge \alpha_i^2\beta_i d_i \overline{\sigma}^2
\end{align}
\end{subequations}
Therefore, we have $\Upsilon_{1}-\Upsilon_{2}^\T\Upsilon_{3}^{-1}\Upsilon_{2}\ge 0$, which completes the proof. $\hfill\blacksquare$

\subsection{Convergence Analysis}

In light of the contraction of the primal-dual gradient updating dynamics, we analyze the convergence of the distributed primal-dual gradient optimization algorithm. The theoretical results are summarized in the following theorem.

\textit{Theorem 2:} Suppose Assumptions 1 and 3 hold. For the distributed primal-dual gradient optimization characterized in $\left(\ref{PDGD2} \right)$, it holds for all $k\ge1$ that
\begin{equation}\label{convergenceAnalysisTheorem}
\left\|\boldsymbol{u}_i^k - \boldsymbol{u}_i^*\right\| \le \frac{2 \mathcal{Q}_i \sqrt{\overline{\zeta} }}{m_J}\tau^k.
\end{equation}
where $0<\mathcal{Q}_i<\infty$ is a scalar depended on the initial and optimal values of $\boldsymbol{\lambda}$.

\textit{Proof.} By using the triangle inequality $\left\|a+b \right\|\le \left\|a \right\|+\left\|b \right\|$, we can obtain
\begin{equation}\label{convergenceAnalysisProof1}
\left\|\boldsymbol{u}_i^k - \boldsymbol{u}_i^*\right\| \le \left\|\boldsymbol{u}_i^k - \left(\boldsymbol{u}_i^k\right)^*\right\|+\left\|\left(\boldsymbol{u}_i^k\right)^* - \boldsymbol{u}_i^*\right\|.
\end{equation}
It is seen that the condition in $\left(\ref{u_i} \right)$ is exact which implies $\left\|\boldsymbol{u}_i^k - \left(\boldsymbol{u}_i^k\right)^*\right\| = 0$. Thereafter, according to the $m_J$-strong convexity in $\left(\ref{assumption1} \right)$, we have
\begin{equation}\label{convergenceAnalysisProof2}
\left\|\left(\boldsymbol{u}_i^k\right)^* - \boldsymbol{u}_i^*\right\|\le\frac{2}{m_J}\left\|\nabla J_i\left(x_i,\left({\boldsymbol{u}}_i^k\right)^*\right)- \nabla J_i\left(x_i,{\boldsymbol{u}}_i^*\right)\right\|
\end{equation}
According to the KKT conditions $\nabla_{{\boldsymbol{u}}_i^k}J_i\left(x_i,\left({\boldsymbol{u}}_i^k\right)^*\right)+G_i^\T\boldsymbol{\lambda}_i^k =0$ and $\nabla_{{\boldsymbol{u}}_i^*}J_i\left(x_i,{\boldsymbol{u}}_i^*\right)+G_i^\T\boldsymbol{\lambda}_i^* =0$, we can obtain
\begin{equation}\label{convergenceAnalysisProof3}
\begin{aligned}
\left\|\left(\boldsymbol{u}_i^k\right)^* - \boldsymbol{u}_i^*\right\|&\le\frac{2}{m_J}\left\|G_i^\T \boldsymbol{\lambda}^k- G_i^\T \boldsymbol{\lambda}^*\right\|\\
&\le\frac{2\sqrt{G_i G_i^\T}}{m_J}\left\|\boldsymbol{\lambda}^k- \boldsymbol{\lambda}^*\right\|\\
&\le\frac{2\mathcal{Q}_i \sqrt{\overline{\zeta}} }{m_J}\tau^k
\end{aligned}
\end{equation}
with $0<\mathcal{Q}_i<\infty$ being a proper scalar depended on the initial and optimal values of $\boldsymbol{\lambda}$, where the last inequality follows from Assumption 3 and the contraction analysis on the PDGD in Theorem 1. The condition in $\left(\ref{convergenceAnalysisTheorem} \right)$ can be established, and the proof is completed. $\hfill\blacksquare$

\subsection{The Stopping Criterion}

In this paper, we adopt the constraint tightening approach to provide the capability of early termination for the DMPC optimization with guaranteeing convergence. The following definition is given to characterize the inexact solution resulted from the early termination. Thereafter, to fulfill the requirement on premature termination, the stopping criterion of the DMPC optimization is developed. 

\textit{Definition 3:} \cite{kohler2019distributed} For a given $\epsilon>0$, if it holds that
\begin{equation}\label{stoppingCriterion1}
\sum_{i=1}^{l} f_i\left(x_i,\boldsymbol{u}_i\right)- b\left(\epsilon \right)\le \epsilon l\boldsymbol{1}_{pN},
\end{equation}
where $\boldsymbol{u}_i\in\tilde{\mathcal{U}}_i\left(x_i\right)$, $\boldsymbol{u}_i$ is called $\epsilon$-feasible solution of the optimization problem in $\left(\ref{erMPC} \right)$.

According to the convergence analysis results on PDGD in $\left(\ref{PDGD2} \right)$, we can establish the stopping criterion in the following theorem.

\textit{Theorem 3:} Given the initial parameters $\epsilon$ and $\gamma_i^0$, the control input sequence $\tilde{u}_i^{\bar{k}}$ is $\epsilon$-feasible under the distributed primal-dual gradient optimization in $\left(\ref{PDGD2} \right)$, where the superscript $\bar{k}$ represents the stopping iteration defined as
\begin{equation}\label{stoppingCriterion}
\bar{k} \ge \log_{\tau} \frac{\epsilon pN}{\left\|\boldsymbol{\gamma}_i^0\right\|}
\end{equation}
with $\tau$ given in Theorem 1, and $\boldsymbol{\gamma}_i^0$ being the initial value of $\boldsymbol{\gamma}_i$.

\textit{Proof.} In light of the PDGD in $\left(\ref{PDGD2} \right)$ and the KKT conditions in $\left(\ref{KKT} \right)$, it holds that
\begin{equation}\label{stoppingCriterion2}
\frac{1}{l}\left(\sum_{i=1}^{l} f_i\left(x_i,\boldsymbol{u}_i\right)- b\left(\epsilon \right)\right) = -\nabla \Psi_i\left(\boldsymbol{\lambda}_i\right)=\boldsymbol{\gamma}_i.
\end{equation}
Therefore, we can rewrite $\left(\ref{stoppingCriterion1} \right)$ in a distributed form as follows
\begin{equation}\label{stoppingCriterion3}
\boldsymbol{\gamma}_i^{\bar{k}} \le \epsilon \boldsymbol{1}_{pN},
\end{equation}
where the superscript $\bar{k}$ stands for the stopping iteration. According to the contraction analysis results in Theorem 1, one can get
\begin{equation}\label{stoppingCriterion4}
\left\|\gamma_i^{\bar{k}}-\gamma_i^* \right\|\le \tau^{\bar{k}}\left\|\gamma_i^0-\gamma_i^* \right\|.
\end{equation}
Following from $\left(\ref{stoppingCriterion3}\right)$ and substituting $\gamma_i^* = 0$ into $\left(\ref{stoppingCriterion4}\right)$, Eq. $\left(\ref{stoppingCriterion}\right)$ can be derived, which completes the proof.$\hfill\blacksquare$

\textit{Remark 5:} It can be seen that the stopping iteration $\overline{k}$ in $\left(\ref{stoppingCriterion}\right)$ depends on the initial values of $\gamma_i$ and the user-defined parameter $\epsilon$. For a given $\epsilon$, to initialize the appropriate values of $\gamma_i$ directly impacts the stopping iteration.

\section{Algorithm Description}

\subsection{Distributed Primal-Dual Gradient Algorithm}

In terms of Theorem 3, the distributed primal-dual gradient algorithm terminates at the iteration $k$, which implies the $\epsilon$-feasible solution is obtained. Thus, the distributed primal-dual gradient algorithm based on Laplacian consensus is described in Algorithm 1. 

\begin{table}[htbp]
	\begin{tabular}{l}
		\toprule
		\textbf{Algorithm 1} The Distributed Primal-Dual Gradient Algorithm \\
		\midrule
		\textbf{Input:} $x_i\left(i=1,\cdots,l\right)$\\
		\textbf{Output:} $\tilde{\boldsymbol{u}}_i\left(i=1,\cdots,l\right)$\\
		\text{Required:} $\epsilon$\\
		\text{Initialization:} $\tilde{\boldsymbol{u}}_i\left(0\right) = 0$, $\boldsymbol{\lambda}_i^0 = 0$ and proper $\boldsymbol{\gamma}_i^0$ for $i=1,\cdots,l$. \\ \quad\quad\quad\quad\quad\quad Obtain stopping criterion $\bar{k}$ based on $\left(\ref{stoppingCriterion4}\right)$;\\
		\;\;1: \textbf{while} $k\le\bar{k}$\\
		\;\;2: \quad \textbf{for} $i=1,\cdots,l$ (in parallel)  \textbf{do}\\
		\;\;3: \quad\quad Exchange $\boldsymbol{\lambda}_i^k$ with its neighbor $j\in\mathcal{N}_i$;\\
		\;\;4: \quad\quad Obtain $\boldsymbol{\gamma}_i^{k+1}$ from $\left(\ref{PDGD2:gamma}\right)$, $\boldsymbol{\lambda}_i^{k+1}$ according to $\left(\ref{PDGD2:lambda}\right)$,\\
		\quad\quad\quad\quad $\tilde{\boldsymbol{u}}_i$ following $\left(\ref{u_i}\right)$, respectively;\\
		\;\;5: \quad \textbf{end for}\\
		\;\;6: \quad Move $k$ to $k+1$.\\
		\;\;7: \textbf{end while}\\
		\bottomrule
	\end{tabular}
\end{table}

\subsection{DMPC Algorithm}

The overall DMPC algorithm adopting the proposed distributed primal-dual gradient algorithm is formulated in the following Algorithm 2.

\begin{table}[htbp]
	\begin{tabular}{l}
		\toprule
		\textbf{Algorithm 2} The DMPC Algorithm \\
		\midrule
		\;\;1: At time instant $t$, each subsystem $i$ samples its state $x_i\left(t\right)$;\\
		\;\;2: Each subsystem $i$ obtains $u_i\left(t\right)$ by following Algorithm 1\\
		\quad\quad with $x_i\left(t\right)$;\\
		\;\;3: Each subsystem $i$ apples its control input;\\
		\;\;4: Let $t=t+1$ and go to step 1 until next sampling instant. \\
		\bottomrule
	\end{tabular}
\end{table}

\section{Recursive Feasibility and Stability}

Under the proposed distributed primal-dual gradient algorithm, the recursive feasibility and stability of the DMPC are analyzed. The following theorem summarized the theoretical analysis results.

\textit{Theorem 4:} Suppose the DMPC optimization problem in $\left(\ref{erMPC} \right)$ is feasible at time instant $t$ for each subsystem $i$ with $i= 1,2,\cdots,l$. Under Assumption 3, the following results hold: i) the DMPC optimization problem in $\left(\ref{erMPC} \right)$ has a feasible solution at time instant $t+1$, ii) the state trajectory of the closed-loop system in $\left(\ref{sysForm} \right)$ enters the terminal set $\Omega_i$ in finite time and remains in it.

\textit{Prrof.} The proof of this theorem consists of two parts known as feasibility and stability analysis, respectively.

\textbf{Feasibility:} Recalling the definition on $\epsilon$-feasible solution in $\left(\ref{stoppingCriterion1} \right)$ and the stopping criterion in $\left(\ref{stoppingCriterion3} \right)$, it is shown that a $\epsilon$-feasible solution $\boldsymbol{u}_i$ satisfies that
\begin{equation}\label{feasiblity1}
\begin{aligned}
\sum_{i=1}^{l}\Phi_i^x x_i\left(s|t\right)+\Phi_i^u u_i\left(s|t\right)&\le \left(1-\epsilon l \left(s+1\right) \right)\boldsymbol{1}_p+\epsilon l\boldsymbol{1}_p \\
&=\left(1-\epsilon l s \right)\boldsymbol{1}_p
\end{aligned}
\end{equation}
with $s\in\mathbb{Z}_{\left[0,N-1\right]}$ according to the coupled constraint in $\left(\ref{coupledConstraint} \right)$ where $\left\lbrace u_i\left(\cdot|t\right)\right\rbrace = \left\lbrace u_i\left(0|t\right),u_i\left(1|t\right),\cdots,u_i\left(N-1|t\right) \right\rbrace$ and $\left\lbrace x_i\left(\cdot|t\right)\right\rbrace = \left\lbrace x_i\left(0|t\right),x_i\left(1|t\right),\cdots,x_i\left(N|t\right) \right\rbrace$. Define a feasible solution at time instant $t+1$ for subsystem $i$ as
\begin{equation}\label{feasiblity2}
\begin{aligned}
&\left\lbrace \tilde{u}_i\left(\cdot|t+1\right)\right\rbrace \\
&:= \left\lbrace \tilde{u}_i\left(0|t+1\right),\tilde{u}_i\left(1|t+1\right),\cdots,\tilde{u}_i\left(N-1|t+1\right) \right\rbrace\\
&:=\left\lbrace u_i\left(1|t\right),u_i\left(2|t\right),\cdots,u_i\left(N-1|t\right), K x_i\left(N|t\right) \right\rbrace,
\end{aligned}
\end{equation}
and the corresponding state sequence can be expressed by
\begin{equation}\label{feasiblity3}
\begin{aligned}
&\left\lbrace \tilde{x}_i\left(\cdot|t+1\right)\right\rbrace \\
&:= \left\lbrace \tilde{x}_i\left(0|t+1\right),\tilde{x}_i\left(1|t+1\right),\cdots,\tilde{x}_i\left(N|t+1\right) \right\rbrace\\
&:=\left\lbrace x_i\left(1|t\right),x_i\left(2|t\right),\cdots,x_i\left(N-1|t\right), A_K x_i\left(N|t\right) \right\rbrace,
\end{aligned}
\end{equation}
where $A_K$ is defined in $\left(\ref{assumption3} \right)$. Thus, we can obtain
\begin{subequations}\label{feasiblity4}
\begin{align}
\sum_{i=1}^{l}\Phi_i^x& \tilde{x}_i\left(s|t+1\right)+\Phi_i^u \tilde{u}_i\left(s|t+1\right) \nonumber\\
&= \sum_{i=1}^{l}\Phi_i^x x_i\left(s+1|t\right)+\Phi_i^u u_i\left(s+1|t\right)\nonumber\\
&\le\left(1-\epsilon l \left(s+1\right) \right)\boldsymbol{1}_p\\
\sum_{i=1}^{l}\left(\Phi_i^x\right.\hspace{-0.3em}&\left. +\Phi_i^u K x_i\left(N|t\right)\right) \le\left(1-\epsilon l N \right)\boldsymbol{1}_p
\end{align}
\end{subequations}
with $s\in\mathbb{Z}_{\left[0,N-2\right]}$. The last inequality follows from the coupled constraint in the terminal set $\left(\ref{localCoupleConstraint} \right)$. Moreover, as the feasibility of the solution at time instant $t$, it holds that $u_i\left(s|t\right)\in\mathcal{U}_i$ and $x_i\left(s|t\right)\in\mathcal{X}_i$ for $s\in\mathbb{Z}_{\left[0,N-1\right]}$ according to the constraint in $\left(\ref{erMPC:U} \right)$. It implies that $\tilde{u}_i\left(s+1|t+1\right)\in\mathcal{U}_i$ and $\tilde{x}_i\left(s+1|t+1\right)\in\mathcal{X}_i$ for $s\in\mathbb{Z}_{\left[0,N-2\right]}$. On the other hand, the local control law $\tilde{u}_i\left(N|t+1\right) = K x_i\left(N|t\right)\in\mathcal{U}_i$ in terms of $x_i\left(N|t\right)\in\Omega_i$. According to the aforementioned analysis results, $\left\lbrace\tilde{u}_i\left(\cdot|t+1\right)\right\rbrace$ is a feasible solution to the DMPC optimization problem in $\left(\ref{erMPC} \right)$ at the successor time instant $t+1$.

\textbf{Stability:} The stability is analyzed by taking the objective function in $\left(\ref{erMPC:ObjFcn}\right)$ as a Lyapunov candidate which refers to $V\left(x_i\left(t\right)\right) := \sum_{i=1}^{l}J_i\left(x_i\left(t\right),\boldsymbol{u}_i^*\left(\cdot|t\right)\right)$. At time instant $t$, denote $\tilde{k}\left(t\right)$ as the iteration $k$ when the stopping criterion in $\left(\ref{stoppingCriterion3} \right)$ is satisfied. In light of Assumption 3, it holds that
\begin{equation}\label{stability1}
\begin{aligned}
J_i&\left(x_i\left(t+1\right),\boldsymbol{u}_i^*\left(\cdot|t+1\right)\right)-J_i\left(x_i\left(t\right),\boldsymbol{u}_i^*\left(\cdot|t\right)\right)\\
&=F\left(x_i\left(N|t\right),Kx_i\left(N|t\right)\right)-F\left(x_i\left(t\right),u_i\left(t\right)\right)\\
&\quad+V_f\left(A_Kx_i\left(N|t\right)\right)-V_f\left(x_i\left(N|t\right)\right)\\
&\le-F\left(x_i\left(t\right),u_i\left(t\right)\right).
\end{aligned}
\end{equation}
As the inexactness of the optimization solver induced by the premature termination, one can have 
\begin{equation}\label{stability2}
\begin{aligned}
V\left(x_i\left(t+1\right)\right) &\le \sum_{i=1}^{l}J_i\left(x_i\left(t+1\right),\boldsymbol{u}_i^*\left(\cdot|t+1\right)\right)\\
& \le \sum_{i=1}^{l}\left(J_i\left(x_i\left(t\right),\boldsymbol{u}_i^*\left(\cdot|t\right)\right)-F\left(x_i\left(t\right),u_i\left(t\right)\right)\right)\\
& = V\left(x_i\left(t\right)\right) -\sum_{i=1}^{l}\left(F\left(x_i\left(t\right),u_i\left(t\right)\right)\right).
\end{aligned}
\end{equation}
It can be seen that $F\left(x_i\left(t\right),u_i\left(t\right)\right)\rightarrow 0$ as $t\rightarrow 0$, which implies there exist finite time instants before the state trajectory enters the terminal set. Thereafter, according to the local controller in Assumption 3, the state trajectory of the closed-loop system will remain in the terminal set, which completes the proof. $\hfill\blacksquare$

\section{Simulation}

\begin{figure}[htbp]
	\centering{\includegraphics[width=8cm,height=2.5cm]{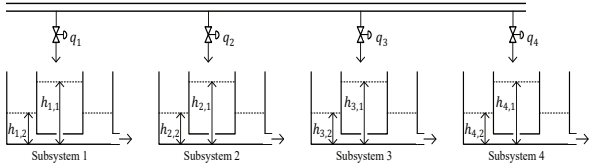}}
	\caption{Double-water tanks systems.}
	\label{simuSysModel}
\end{figure}
\begin{figure}[htbp]
	\centering{\includegraphics[width=8cm,height=4cm]{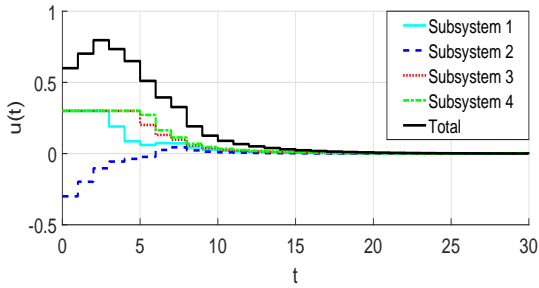}}
	\caption{Control input of each subsystem and the total control input.}
	\label{simu:input}
\end{figure}
\begin{figure}[htbp]
	\centering{\includegraphics[width=8cm,height=6cm]{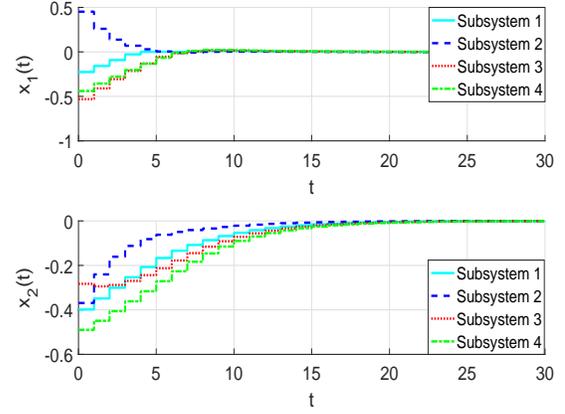}}
	\caption{States of each subsystem.}
	\label{simu:state}
\end{figure}

A numerical simulation is given to verify the performance of the presented DMPC approach. Consider a coupled double-water tank systems containing four subsystems \cite{wang2018accelerated,wellstead1979introduction}. The control objective is to regulate the water level towards the
given reference by means of the input flow which is subject to a globally coupled rate constraint. As shown in Fig.$\ref{simuSysModel}$, denote the sampled water level as $h_{i,1}$ and $h_{i,2}$, the input flow $q_{i}$ for each subsystem $i$ with $i=1,2,3,4$. Given the reference water level as $r_{i,1}=1$ and $r_{i,1}=0.64$, and the steady-state input flow $\overline{q}_{i}=0.3$, formulate each subsystem $i$ as the following discrete-time linear dynamical system
\begin{equation}\label{SimuModel}
x_i\left(t+1\right) = A x_i\left(t\right)+B u_i\left(t\right),
\end{equation}
where $x_{i,1} = h_{i,1}-r_{i,1}$, $x_{i,2} = h_{i,2}-r_{i,2}$, $u_{i} = q_{i}-\overline{q}_{i}$, and
\begin{equation}\label{SimuModel_AB}
A = 
\begin{bmatrix}
0.8750&0.1250\\
0.1250&0.8047
\end{bmatrix},
B = 
\begin{bmatrix}
0.3\\
0
\end{bmatrix}.
\end{equation}
Define the objective function for the DMPC optimization problem as
\begin{figure}[htbp]
	\centering{\includegraphics[width=8cm,height=5cm]{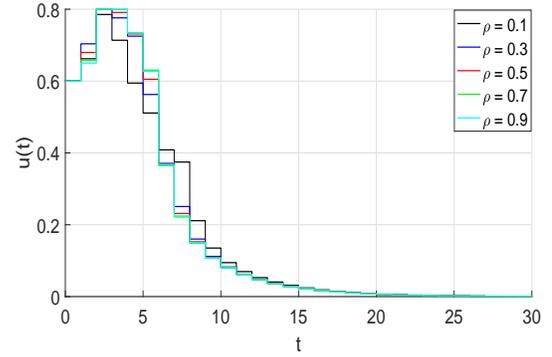}}
	\caption{Comparison on $\rho$.}
	\label{simu:Compare_rho}
\end{figure}
\begin{figure}[htbp]
	\centering{\includegraphics[width=8cm,height=5cm]{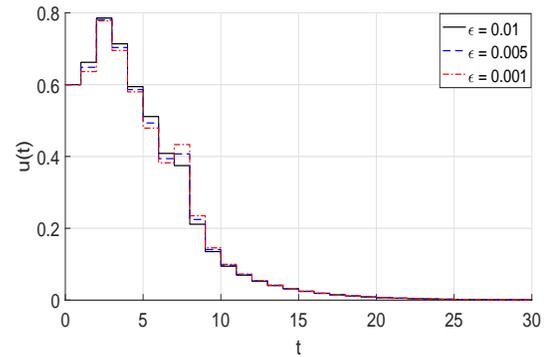}}
	\caption{Comparison on $\epsilon$.}
	\label{simu:Compare_epsilon}
\end{figure}
\begin{equation}\label{SimuObjFcn}
J_i\left(x_i,u_i\right) = \sum_{s=0}^{N}\left( \left\|x_i\left(s|t\right)\right\|_{Q_i}^2+\left\|u_i\left(s|t\right)\right\|_{R_i}^2\right)+\left\|x_i\left(N|t\right)\right\|_{P_i}^2
\end{equation}
where $N = 8$, $Q_i = 10\boldsymbol{I}$, $R_i = 1$,
\begin{equation}\label{PQR}
P_i = 
\begin{bmatrix}
31.7459&9.8300\\
9.8300&56.3415
\end{bmatrix}
\end{equation}
and the local control gain matrix $K_i = \left[-1.7916\,-0.7337 \right]$. Taking the physical constraints into account, each subsystem should consider the local constraints on states as $\mathcal{X}_i = \left\lbrace x_i\in\mathbb{R}^2|-1\le x_{i,1} \le 1, -0.64 \le x_{i,2} \le 0.64\right\rbrace$ and input $\mathcal{U}_i = \left\lbrace u_i\in\mathbb{R}|-0.3\le u_{i} \le 0.3\right\rbrace$, respectively. As the presence of the bound on total input flow rate which is supposed to be 2 in this paper, a global constraint is involved in the DMPC optimization problem as $\sum_{i=1}^{4}q_{i}\le 2$, which means $\sum_{i=1}^{4}\Phi_i^u u_i \le 1$ with $\Phi_i^u = 1.25$. In addition, the adjacent matrix and the corresponding Laplacian matrix of the connection network are given as follows
\begin{equation}\label{SimuAdjacentMatrix}
\mathcal{A} = 
\begin{bmatrix}
0&1&0&1\\
1&0&1&0\\
0&1&0&1\\
1&0&1&0
\end{bmatrix},
L = 
\begin{bmatrix}
2&-1&0&-1\\
-1&2&-1&0\\
0&-1&2&-1\\
-1&0&-1&2
\end{bmatrix}.
\end{equation}
The initial states of each subsystem are $x_1\left(0\right) = [-0.2264\,-0.3981]^\T$, $x_2\left(0\right) = [0.4520\,-0.3689]^\T$, $x_3\left(0\right) = [-0.5311\,-0.2828]^\T$ and $x_4\left(0\right) = [-0.4397\,-0.4897]^\T$, respectively. Define the accuracy for coupled constraints as $\epsilon = 0.1$, and $\boldsymbol{\gamma}_i^0$ is initialized accordingly. In terms of Theorem 1, we choose $\alpha = 0.2$ and $\beta = 0.19$. 

The simulation results are shown in Fig. \ref{simu:input}-Fig. \ref{simu:Compare_epsilon}. Fig. \ref{simu:input} gives the control input of each subsystems. Meanwhile, the total control input demonstrates the satisfaction of the global coupled constraints. In Fig. \ref{simu:state}, the states of each subsystem are plotted, which shows that the control objective can be achieved. To illustrate the effectiveness of the user-defined parameters introduced in the proposed approach, we give the total control input under different $\rho$ and $\epsilon$ in Fig. \ref{simu:Compare_rho} and Fig. \ref{simu:Compare_epsilon}, respectively. It is seen that the changing rate of the control input is increasing as $\rho$ gets bigger. Moreover, a smaller $\epsilon$ leads to more conservative of control input.

\section{Conclusion}

In this paper, a novel distributed model predictive control (DMPC) approach is studied for a group of discrete-time linear systems taken the global coupled constraints into account. The DMPC optimization problem is transformed into a dual problem involving all subsystems, which is solved in the framework of Laplacian consensus by using the primal-dual gradient optimization in a fully distributed manner. To reduce the computational burden, a tightening constraint concerning the global coupled constraint is constructed to allow premature termination with guaranteeing the convergence of the optimization. It is worth seeing that the local copies of the Lagrangian multipliers need not reach consensus but within some specified bounds owing to the constraint tightening method. Furthermore, the convergence of the primal-dual gradient optimization is first rigorously analyzed by means of contraction theory in the context of discrete-time nonlinear dynamics. The recursive feasibility and stability of the closed-loop system are established under the inexact solver with rational assumptions. The performance of the proposed approach is demonstrated by a numerical simulation.

% or


%\appendices

%\section{}

%\section*{Acknowledgment}

\ifCLASSOPTIONcaptionsoff
  \newpage
\fi



%\bibliographystyle{IEEEtran}
\bibliography{mybibfile}

\begin{thebibliography}{10}
\providecommand{\url}[1]{#1}
\csname url@samestyle\endcsname
\providecommand{\newblock}{\relax}
\providecommand{\bibinfo}[2]{#2}
\providecommand{\BIBentrySTDinterwordspacing}{\spaceskip=0pt\relax}
\providecommand{\BIBentryALTinterwordstretchfactor}{4}
\providecommand{\BIBentryALTinterwordspacing}{\spaceskip=\fontdimen2\font plus
\BIBentryALTinterwordstretchfactor\fontdimen3\font minus
  \fontdimen4\font\relax}
\providecommand{\BIBforeignlanguage}[2]{{%
\expandafter\ifx\csname l@#1\endcsname\relax
\typeout{** WARNING: IEEEtran.bst: No hyphenation pattern has been}%
\typeout{** loaded for the language `#1'. Using the pattern for}%
\typeout{** the default language instead.}%
\else
\language=\csname l@#1\endcsname
\fi
#2}}
\providecommand{\BIBdecl}{\relax}
\BIBdecl

\bibitem{qin2003survey}
S.~J. Qin and T.~A. Badgwell, ``A survey of industrial model predictive control
  technology,'' \emph{Control Engineering Practice}, vol.~11, no.~7, pp.
  733--764, 2003.

\bibitem{eren2017model}
U.~Eren, A.~Prach, B.~B. Ko{\c{c}}er, S.~V. Rakovi{\'c}, E.~Kayacan, and
  B.~A{\c{c}}{\i}kme{\c{s}}e, ``Model predictive control in aerospace systems:
  Current state and opportunities,'' \emph{Journal of Guidance, Control, and
  Dynamics}, vol.~40, no.~7, pp. 1541--1566, 2017.

\bibitem{liu2018robust}
C.~Liu, H.~Li, J.~Gao, and D.~Xu, ``Robust self-triggered min--max model
  predictive control for discrete-time nonlinear systems,'' \emph{Automatica},
  vol.~89, pp. 333--339, 2018.

\bibitem{kvasnica2019complexity}
M.~Kvasnica, P.~Bakar{\'a}{\v{c}}, and M.~Klau{\v{c}}o, ``Complexity reduction
  in explicit {MPC}: A reachability approach,'' \emph{Systems \& Control
  Letters}, vol. 124, pp. 19--26, 2019.

\bibitem{su2018self}
Y.~Su, Q.~Wang, and C.~Sun, ``Self-triggered robust model predictive control
  for nonlinear systems with bounded disturbances,'' \emph{IET Control Theory
  \& Applications}, vol.~13, no.~9, pp. 1336--1343, 2019.

\bibitem{li2017robust}
H.~Li and Y.~Shi, \emph{Robust Receding Horizon Control for Networked and
  Distributed Nonlinear Systems}.\hskip 1em plus 0.5em minus 0.4em\relax New
  York, NY, USA: Springer, 2017.

\bibitem{li2016neighbor}
H.~Li, Y.~Shi, and W.~Yan, ``On neighbor information utilization in distributed
  receding horizon control for consensus-seeking,'' \emph{IEEE Transactions on
  Cybernetics}, vol.~46, no.~9, pp. 2019--2027, 2016.

\bibitem{chen2018adaptive}
Y.~Chen, M.~Bruschetta, D.~Cuccato, and A.~Beghi, ``An adaptive partial
  sensitivity updating scheme for fast nonlinear model predictive control,''
  \emph{IEEE Transactions on Automatic Control}, in press, 2018.

\bibitem{kohler2019distributed}
J.~K{\"o}hler, M.~A. M{\"u}ller, and F.~Allg{\"o}wer, ``Distributed model
  predictive control—recursive feasibility under inexact dual optimization,''
  \emph{Automatica}, vol. 102, pp. 1--9, 2019.

\bibitem{dunbar2012distributed}
W.~B. Dunbar and D.~S. Caveney, ``Distributed receding horizon control of
  vehicle platoons: Stability and string stability,'' \emph{IEEE Transactions
  on Automatic Control}, vol.~57, no.~3, pp. 620--633, 2012.

\bibitem{kohler2018distributed}
P.~N. K{\"o}hler, M.~A. M{\"u}ller, and F.~Allg{\"o}wer, ``A distributed
  economic {MPC} framework for cooperative control under conflicting
  objectives,'' \emph{Automatica}, vol.~96, pp. 368--379, 2018.

\bibitem{liu2018distributed}
P.~Liu and U.~Ozguner, ``Distributed model predictive control of spatially
  interconnected systems using switched cost functions,'' \emph{IEEE
  Transactions on Automatic Control}, vol.~63, no.~7, pp. 2161--2167, 2018.

\bibitem{shi2015extra}
W.~Shi, Q.~Ling, G.~Wu, and W.~Yin, ``{EXTRA}: An exact first-order algorithm
  for decentralized consensus optimization,'' \emph{SIAM Journal on
  Optimization}, vol.~25, no.~2, pp. 944--966, 2015.

\bibitem{fazlyab2018distributed}
M.~Fazlyab, S.~Paternain, A.~Ribeiro, and V.~M. Preciado, ``Distributed smooth
  and strongly convex optimization with inexact dual methods,'' in \emph{2018
  Annual American Control Conference (ACC)}.\hskip 1em plus 0.5em minus
  0.4em\relax IEEE, 2018, pp. 3768--3773.

\bibitem{patrascu2018convergence}
A.~Patrascu and I.~Necoara, ``On the convergence of inexact projection primal
  first-order methods for convex minimization,'' \emph{IEEE Transactions on
  Automatic Control}, vol.~63, no.~10, pp. 3317--3329, 2018.

\bibitem{yu2018convergence}
H.~Yu and M.~J. Neely, ``On the convergence time of dual subgradient methods
  for strongly convex programs,'' \emph{IEEE Transactions on Automatic
  Control}, vol.~63, no.~4, pp. 1105--1112, 2018.

\bibitem{hale2017asynchronous}
M.~T. Hale, A.~Nedi{\'c}, and M.~Egerstedt, ``Asynchronous multiagent
  primal-dual optimization,'' \emph{IEEE Transactions on Automatic Control},
  vol.~62, no.~9, pp. 4421--4435, 2017.

\bibitem{dunbar2007distributed}
W.~B. Dunbar, ``Distributed receding horizon control of dynamically coupled
  nonlinear systems,'' \emph{IEEE Transactions on Automatic Control}, vol.~52,
  no.~7, pp. 1249--1263, 2007.

\bibitem{dai2017distributed}
L.~Dai, Y.~Xia, Y.~Gao, and M.~Cannon, ``Distributed stochastic {MPC} of linear
  systems with additive uncertainty and coupled probabilistic constraints,''
  \emph{IEEE Transactions on Automatic Control}, vol.~62, no.~7, pp.
  3474--3481, 2017.

\bibitem{wang2017distributed}
Z.~Wang and C.~J. Ong, ``Distributed model predictive control of linear
  discrete-time systems with local and global constraints,'' \emph{Automatica},
  vol.~81, pp. 184--195, 2017.

\bibitem{Lohmiller1998}
W.~Lohmiller and J.~J.~E. Slotine, ``{On contraction analysis for non-linear
  systems},'' \emph{Automatica}, vol.~34, no.~6, pp. 683--696, 1998.

\bibitem{long2018distributed}
Y.~Long, S.~Liu, L.~Xie, and K.~H. Johansson, ``Distributed nonlinear model
  predictive control based on contraction theory,'' \emph{International Journal
  of Robust and Nonlinear Control}, vol.~28, no.~2, pp. 492--503, 2018.

\bibitem{pham2009contraction}
Q.~C. Pham, N.~Tabareau, and J.~J. Slotine, ``A contraction theory approach to
  stochastic incremental stability,'' \emph{IEEE Transactions on Automatic
  Control}, vol.~54, no.~4, pp. 816--820, 2009.

\bibitem{forni2014differential}
F.~Forni and R.~Sepulchre, ``A differential {Lyapunov} framework for
  contraction analysis.'' \emph{IEEE Transactions on Automatic Control},
  vol.~59, no.~3, pp. 614--628, 2014.

\bibitem{chaffey2018control}
T.~L. Chaffey and I.~R. Manchester, ``Control contraction metrics on {Finsler}
  manifolds,'' \emph{arXiv preprint arXiv:1803.01034}, 2018.

\bibitem{nguyen2018contraction}
H.~D. Nguyen, T.~L. Vu, K.~Turitsyn, and J.-J. Slotine, ``Contraction and
  robustness of continuous time primal-dual dynamics,'' \emph{IEEE Control
  Systems Letters}, vol.~2, no.~4, pp. 755--760, 2018.

\bibitem{liu2017robust}
X.~Liu, Y.~Shi, and D.~Constantinescu, ``Robust distributed model predictive
  control of constrained dynamically decoupled nonlinear systems: A contraction
  theory perspective,'' \emph{Systems \& Control Letters}, vol. 105, pp.
  84--91, 2017.

\bibitem{stuart1994numerical}
A.~M. Stuart, ``Numerical analysis of dynamical systems,'' \emph{Acta
  Numerica}, vol.~3, pp. 467--572, 1994.

\bibitem{qu2018exponential}
G.~Qu and N.~Li, ``On the exponential stability of primal-dual gradient
  dynamics,'' \emph{IEEE Control Systems Letters}, vol.~3, no.~1, pp. 43--48,
  2018.

\bibitem{rubagotti2014stabilizing}
M.~Rubagotti, P.~Patrinos, and A.~Bemporad, ``Stabilizing linear model
  predictive control under inexact numerical optimization,'' \emph{IEEE
  Transactions on Automatic Control}, vol.~59, no.~6, pp. 1660--1666, 2014.

\bibitem{wang2018accelerated}
Z.~Wang and C.~J. Ong, ``Accelerated distributed {MPC} of linear discrete-time
  systems with coupled constraints,'' \emph{IEEE Transactions on Automatic
  Control}, vol.~63, no.~11, pp. 3838--3849, 2018.

\bibitem{hashimoto2016self}
K.~Hashimoto, S.~Adachi, and D.~V. Dimarogonas, ``Self-triggered model
  predictive control for nonlinear input-affine dynamical systems via adaptive
  control samples selection,'' \emph{IEEE Transactions on Automatic Control},
  vol.~62, no.~1, pp. 177--189, 2016.

\bibitem{bertsekas1999nonlinear}
D.~P. Bertsekas, \emph{Nonlinear Programming}.\hskip 1em plus 0.5em minus
  0.4em\relax Belmont, MA, USA: Athena Scientific, 1999.

\bibitem{wellstead1979introduction}
P.~E. Wellstead, \emph{Introduction to Physical System Modelling}.\hskip 1em
  plus 0.5em minus 0.4em\relax Cambridge, MA, USA: Academic Press, 1979.

\end{thebibliography}


%\end{IEEEbiography}
%% insert where needed to balance the two columns on the last page with
%%\newpage


%\vfill


\end{document}